\documentclass[12pt]{article}
\usepackage{amsmath}
\usepackage{amssymb}
\usepackage{latexsym}
\usepackage{epsfig}
\usepackage{graphicx}
\usepackage{mathptmx}
\usepackage[mtpcal]{}
\usepackage{mathtools}
\usepackage{physics}
\usepackage{cancel}
\usepackage{appendix}

\begin{document}
	\setlength{\baselineskip}{7mm}

	\setlength{\baselineskip}{10mm}
	
	\begin{center}
		{\huge Deriving Lehmer and H\"older means as maximum
			weighted likelihood estimates for the multivariate
			exponential family}         
	\end{center}
	\begin{center}
			Djemel Ziou\footnote{Corresponding author:	Djemel.Ziou@usherbrooke.ca} and Issam Fakir
			
			\vspace*{-4mm}
			D\'epartement d'informatique  
			
			\vspace*{-4mm}
			Universit\'e de Sherbrooke 
			
			\vspace*{-4mm}
			Sherbrooke, Qc., Canada J1K 2R1
	\end{center}
	
	\setlength{\baselineskip}{7mm}
	
	
	\vspace*{5mm}

\begin{abstract}
The links between the mean families of Lehmer and H\"older and the weighted maximum likelihood estimator have recently been established in the case of a regular univariate exponential family. In this article, we will extend the outcomes obtained to the multivariate case. This extension provides a probabilistic interpretation of these families of means and could therefore broaden their uses in various applications.
\end{abstract}


\section{Introduction}
Consider numerical observations; it is common to calculate their mean and refer to it as central tendency. There are, however, different measures of mean \cite{Beckenbach50}. These measurements are sometimes grouped into families, like Lehmer and H\"older. Distinguishing these measures  and better understanding their use involves identifying the link between them and probability density functions (PDFs). For example, the arithmetic mean is the maximum likelihood estimator (MLE) of the position parameter for the normal PDF and the scale parameter for the exponential PDF. For the families of Lehmer and H\"older means, such an interpretation has only recently been proposed for the case of PDFs in the case of the univariate exponential family Let’s consider digital observations; it is often common to calculate their mean and designate it as a central tendency. However, there are various measures of the average \cite{mean_ge}. These measures are sometimes grouped into families, such as Lehmer and H\"older. Distinguishing these mean measures and gaining a deeper understanding of their use involves identifying the connection between them and the probability density functions (PDFs). For instance, the arithmetic mean is the maximum likelihood estimator of the position parameter for the normal PDF and the scale parameter for the exponential PDF. For the Lehmer and H\"older families of means, such an interpretation has only recently been proposed for the case of PDFs in the univariate exponential family case \cite{Ziou23a}.
\newline
In the following text, the terms "mean", "mean family", and “central tendency” are used interchangeably. This paper extends the work in~\cite{Duerinckx14} to the multivariate case. We establish the relationship between the weighted maximum likelihood estimator (MWLE) and the Lehmer and H\"older family of means in the context of the multivariate exponential family. Specifically, we demonstrate that: 1) these means are MWLEs for a subset of this family of PDFs; 2) the MWLE not only depends on a PDF, as shown in existing maximum likelihood estimator characterization studies~\cite{Ziou23a}, but also on the relevance of the data. These correspondences offer a probabilistic explanation for these means, potentially broadening their applicability in different domains that rely on maximum likelihood estimation. The paper is organized as follows. The next section outlines the two central tendencies. In section \ref{MaximumLikelihhodEsi}, we derive MWLEs and illustrate their connection to these central tendencies. Section \ref{case studie} presents a case study.

\section{Lehmer's and H\"older's Central tendencies}
\begin{figure}
	\begin{center}
		\includegraphics[height=4cm, width=4cm]{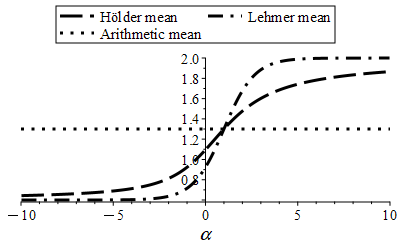}{a)}
		\includegraphics[height=4cm, width=4cm]{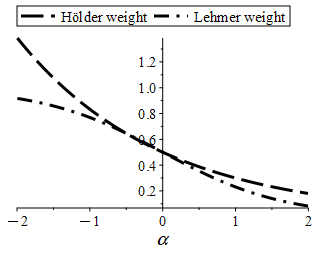}{b)}
		\includegraphics[height=4cm, width=4cm]{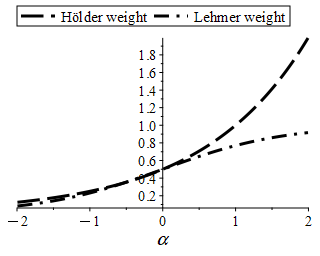}{c)}
		\caption{(a) Three means between $0.6$ and $2$ as function of $\alpha$. The arithmetic mean is displayed as a basis for comparison. (b) The weights  $m_h(0.6)$ and $m_l(0.6)$  as a function of $\alpha$.  (c) The weights  $m_h(2)$ and $m_l(2)$ as a function of $\alpha$.   
		}
		\label{MeansPlot}
	\end{center}
\end{figure}

Let us define the notion of centrality for one random variable $X_j$ and its observations ($x_{1,j},...,x_{n,j}$),
where $(x_{i,j})_{1\leq i\leq n} \in  \mathbb{R}_{\ge 0}$. Several formulas exist for calculating the mean of our observations, and most of them are  a particular case of the generalized f-mean, known also as the generalized Kolmogorov mean~\cite{generalization}:
\begin{equation}
	 \mu_j = f^{-1} (\frac{1}{n}\sum_{i=1}^n f(x_{i,j}))
	\label{Kolmogorov}
\end{equation}
where the function $f : \mathbb{R}_{\ge 0} \rightarrow \mathbb{R}_{\ge 0}$  is continuous and increasing. The  H\"older  family of means $H_{\alpha}^{(j)}$ is  particular cases of the generalized f-mean, where $\alpha \in \mathbb{R}$ is the H\"older parameter. It is obtained by setting  $f(x)=x^\alpha$: 
\begin{equation}
	H_{\alpha}^{(j)}=(\frac{\sum_{i=1}^n w(x_{i,j}) x_{i,j}^\alpha}{\sum_{i=1}^n w(x_{i,j})})^{1/\alpha}
	\label{Holder}	
\end{equation}
Where $w$ is weight function. Due to its connection with the $p$-norm, it is widely used in information technology, finance, health,  and human development assessment, pattern recognition, among many other areas~\cite{Tripathi11,deCarvalho16,Oh16}.  Lehmer $L_{\alpha}^{(j)}$ is an alternative to the f-mean and it is  given by: 
\begin{equation}
	L_{\alpha}^{(j)}=\frac{\sum_{i=1}^n w(x_{i,j}) x_{i,j}^\alpha}{\sum_{i=1}^n w(x_{i,j}) x_{i,j}^{\alpha-1}}
	\label{Lehmer}	
\end{equation}
It is used in  differential evolution~\cite{Das17}, neural networks~\cite{Terziyan22},  extreme events estimation~\cite{Penalva20}, and depressive disorders characterization~\cite{Ataei22}.
Both families are bounded by the smallest and greatest values of the sample and  are  continuously non-decreasing  functions with respect to $\alpha$. 
For the purpose of comparison  between both, Fig. \ref{MeansPlot}.a represents  the Lehmer and H\"older  central tendencies of the two equiprobably numbers $x_{1,j}=0.6$ and $x_{2,j}=2$ as function of $\alpha$ and the arithmetic mean as a basis. The Lehmer  is greater than the  H\"older  when $\alpha > 1$,  lower for $\alpha < 1$, and equal when $\alpha=-\infty, 1, +\infty$.  The Pythagorean central tendencies are particular cases;  the geometrical mean $H_0^{(j)} = L_{0.5}^{(j)}$ when $n=2$,  the arithmetic mean $H_1^{(j)}=L_1^{(j)}$, and the harmonic mean 
$H_{-1}^{(j)}=L_0^{(j)}$.  Moreover, the slope of the Lehmer  is higher than the H\"older, i.e. the Lehmer  reaches  the lowest and highest values more quickly. Both are smaller than the arithmetic mean when $\alpha<1$,
greater when $\alpha > 1$ and equal when $\alpha=1$. The reader can find more details about these means in~\cite{Beckenbach50,Burrows86}.
An important issue we want to address is the explanation of the data selection embedded  in the central tendencies.   
Drawing inspiration from the interpretation of the mean described in~\cite{Ziou22a}, we provide two data selection mechanisms  embedded  in each of  Lehmer and H\"older means. The first is the w-weight $w()$ encoding   knowledge about the observations, such as the frequency  or prior knowledge. The higher $w(x_{i,j})$, the more $x_{i,j}$ contributes to the central tendency.  The second mechanism is the v-weight $v(x_{i,j})$ based on the value $x_{i,j}$; that is to say that a measurement contributes to the calculation of the central tendency according to its value $x_{i,j}$. More precisely, the H\"older and Lehmer v-weighted means can be written, respectively, as:  
\begin{equation}
	H_{\alpha}^{(j)}=(\sum_{i=1}^n v_h(x_{i,j}) x_{i,j})^{1/\alpha}~~~\mbox{and}~~~L_{\alpha}^{(j)}=\sum_{i=1}^n v_l(x_{i,j}) x_{i,j}
	\label{Arithmean}
\end{equation}
where  $v_l(x_{i,j})=w(x_{i,j}) x_{i,j}^{\alpha-1}/\sum_{i=1}^n w(x_{i,j}) x_{i,j}^{\alpha-1}$ and $v_h(x_{i,j})=w(x_{i,j}) x_{i,j}^{\alpha-1}/\sum_{i=1}^n w(x_{i,j})$. Note that while the sum of the v-weights is equal to one (i.e. $ \sum_{i=1}^n v_l(x_{i,j})=1$) in the case of Lehmer, it is not the same for H\"older because $\sum_{i=1}^n v_h(x_{i,j})=(H^{(j)}_{\alpha-1})^{\alpha-1}$ for $\alpha \ne 1$. The link between the two means is  straightforward  $L_{\alpha}^{(j)}=(H_{\alpha}^{(j)})^\alpha/(H^{(j)}_{\alpha-1})^{\alpha-1}$.  To better illustrate these v-weights, let us consider again the two equiprobably numbers $x_{1,j}=0.6$ and $x_{2,j}=2$.   Fig. \ref{MeansPlot}.b depicts  $v_h(0.6)=0.6^\alpha/2$ and $v_l(0.6)=0.6^\alpha/ ( 0.6^\alpha+2^\alpha)$    and Fig. \ref{MeansPlot}.c  $v_h(2)= 2^\alpha/2$ and $v_l(2)= 2^\alpha/ ( 0.6^\alpha+2^\alpha) $ as function of $\alpha$. The function $v_z(0.6)$ is decreasing and $v_z(2)$ is increasing in $\alpha$, where $z$ is either $l$ or  $h$. In other words, the relevance of $x_{1,j}=0.6$  (resp. $x_{2,j}=2$)   is decreasing (resp. increasing) when $\alpha$ is increasing. Hence, one of the data selection mechanisms embedded in the two means involves increasing values above one and weakening values below one.

\section{Maximum likelihood estimates\label{MaximumLikelihhodEsi}}
Let us go considering the weighted data \begin{equation}
    {\cal D}= {\cal X} \times {\cal U} = \begin{pmatrix}
x_{i,1}&...&x_{1,k}&u(x_{i,1},...,x_{1,k}) \\ \vdots & \vdots & \vdots & \vdots \\ x_{n,1}&...&x_{n,k}&u(x_{n,1},...,x_{n,k}))\end{pmatrix} \label{Dataset weight}
\end{equation} where $\forall i \in \{1,...,n\}, \ u(x_i) \in \mathbb{R}_{> 0}$ is a weight of the observation $x_i$.  The relationship between the weight function $u(x)$ and the previously introduced w-weight and v-weight will be further elaborated.
Let us consider that the data $\cal H$ as a realization of random vector X of distribution $\mathbb{P}(\theta)$, with $\theta$ in $\Theta$ the space parameter, where the PDF is in the multivariate exponential family supposed minimal. The latter is written as follows:
\begin{equation}
	f(x|\theta)=a(x) exp(<\eta(\theta), T(x)>-H(\theta)) 
	\label{pdfgexp}
\end{equation}
where  $H(\theta)=ln \int_{\mathbb {R}^k} a(x) exp(<\eta(\theta),T(x)>) dx$ is the normalizer, $a(x)$ the basis measure  associating non-negative values to $x$ regardless of $\theta$,  $T(x)\in\mathbb{R}^q$ is referred to as  a sufficient statistic, and $\eta(\theta)\in \mathbb{R}^q$  the parametrization function which is supposed bijective on $\Theta$. The quantity $q$ being the degree of liberty of our random vector, such as the elements of $T(X)$ aren't redundant (see Appendix \ref{minimal_expo}). By setting  each of $a(x)$, $T(x)$, and $\eta(\theta)$ to  a specific value, several existing PDFs can be derived from Eq. \ref{pdfgexp}  such as the multinomial, multivariate Gaussian.
\newline
 We will use the weighted likelihood to establish the relationship between Lehmer and H\"older's central tendencies with the MWLE. The weighted likelihood was proposed by Feifang Hu and which consists of integrating observation weights into the Fisher likelihood~\cite{Hu94}. Under the IID assumption, for the PDFs in Eq. \ref{pdfgexp}, the weighted likelihood of the weighed data ${\cal D}$ relatively to the naturals parameters $\eta \in \{\eta, H(\eta)<\infty\}$ is given by:
\begin{equation}
\begin{split}
    L(\eta)= \prod_i (a(x_i) exp(<\eta,T(x_i)>  - H(\eta)))^{u(x_i)}
\end{split}
\label{LH}
\end{equation}
The positive weight function $u(x)$  must be $\theta$ free. We define the log-weighted likelihood function such as:
\begin{equation}
 \begin{split}
     l_n(\eta)=\sum_{i=1}^n u(x_i)[\ln(a(x_i))+<\eta,T(x_i)> - H(\eta))]
 \end{split}
 \label{ln eta}
\end{equation}
Equating the first derivative with respect to $\eta$ of $l_n()$ to zero and resolving, leads us to write (Appendix \ref{calcul_det}):
\begin{equation}
	r(\eta)=\frac{\sum_{i=1}^n u(x_i) T(x_i)}{\sum_{i=1}^n u(x_i)}
	\label{CP}
\end{equation}
where $r(\eta)=\grad H(\eta)$.  Since the derivative of $H(\eta)$ wrt $\eta$ is $\grad H(\eta)=  \mathbb{E}_{X\sim\eta}[T(X)]$, then   $r(\eta) =\mathbb{E}_{X\sim\eta}[T(X)]$, where $\mathbb{E}_{X\sim\eta}[T(X)]$ is the expectation of $T(X)$ relatively to the PDF of parameter $\eta$ (Eq. \ref{pdfgexp}) and with X a random vector following the distribution $\mathbb{P}(\eta)$. 
By using  the Lebesgue dominated convergence theorem \cite{LDCT}, the  Hessian matrix of $l_n(\eta)$ wrt $\eta$ is given by:
\begin{equation}
	Hess_{l_n}(\eta) =  -(\sum_{i=1}^n u(x_i)) K_{T(X),T(X)}
\end{equation}
Where $K_{T(X),T(X)}$ is the covariance matrix of $T(X)$ relatively to the PDF of parameter $\eta$ which is a positive semi-definite matrix  and $(\sum_{i=1}^n u(x_i))$ the total weight which is positive and non-null. So $l_n()$ is concave, and a critical point of $l_n()$ is also a global maximum of $l_n()$ (which is the MWLE by definition). This global maximum is unique when the density belongs to the minimal exponential family \cite{Expo invertible}. It follows that there is a one-to-one mapping between $r()$ and its global maximum (where the inverse of a vector is the inverse of its components), and, therefore,  the critical point is: 
\begin{equation}
 \hat{\eta} = r^{-1}(\frac{\sum_{i=1}^n u(x_i) T(x_i)}{\sum_{i=1}^n u(x_i))})
\label{CriticalP}
\end{equation}
$\eta$ is supposed bijective so we can find $\hat{\theta}$ by numerical methods or find its form when it's possible, 
\begin{equation}
   \begin{split}
       \hat{\theta} &= \eta^{-1}(\hat{\eta})
       \\&= \eta^{-1}(r^{-1}(\frac{\sum_{i=1}^n u(x_i) T(x_i)}{\sum_{i=1}^n u(x_i)}) )
   \end{split} 
\end{equation}
Unlike Eq. \ref{Kolmogorov}, four functions are involved in this formula that are $\eta()$, $r()$, $u()$, and $T()$. This formula, less specialized than the generalized f-mean in Eq. \ref{Kolmogorov}, will  be used later to derive the Lehmer and H\"older means. 
\label{Lehmer/Holder}
Another important issue concerns the selection of data, namely the function $u(x)$ to use.  In our case, the objective is to use the functions $u(x)$ which make it possible to derive   Lehmer and H\"older means as  MWLE. 
\newline
Comparing Eq. \ref{CP} and H\"older mean in Eq. \ref{Arithmean}  raised to $\alpha$ power gives $u(x_i)=w(x_i)$ and $T(x_i)=(x_{i,1}^{\alpha_1},...,x_{i,q}^{\alpha_q})$ with $k\geqslant q$ (We explain why we need that $k\geqslant q$ further); that is
the MWLE estimator in  Eq. \ref{CP} is a function of H\"older mean. The subclass of PDFs from the multivariate exponential family leading to a function of H\"older mean as MWLE has the form $a(x) exp(\sum_{i=1}^q x_i^{\alpha_i} \eta_i (\theta) -H(\theta))$ when $u(x)= w(x)$.
\newline \label{Holder_Lehmer_u}
For Lehmer's estimator, the random variables need to be independents, comparing Eq. \ref{CP} and  Lehmer mean in Eq. \ref{Arithmean} gives this form $u(x_i)= w(x_{i,j}) x_{i,j}^{\alpha_j-1}$  and $T(x_i)= (x_{i,1},...,x_{i,q})$ with $k\geqslant q$; that is  the MWLE estimator  in Eq. \ref{CP} is a function of Lehmer mean. The subclass of PDFs from the multivariate exponential family leading to a function of Lehmer mean as MWLE has the form $a(x) exp(\sum_{j=1}^q\eta(\theta)x_j-H(\theta))$ when $u(x_{i,j})= w(x_{i,j}) x_{i,j}^{\alpha_j-1}$.
\newline
Let us explain why $k\geq q$: Consider the sufficient statistics and natural parameter of a multinomial distribution Multi(k,N,p) with k and N fixed, represented as:
\begin{equation*}
   \begin{split}
       &\mathbb{P}(X=x_i)=\frac{N!}{x_{i,1}!...x_{i,k}!}exp(\sum_{j=1}^k x_{i,j} ln(p_j) ), \\&\ x_i=(x_{i,j})_{1\leq j\leq k}\in\mathbb{R} \ \text{such as} \ \sum_j x_j=N
   \end{split} 
\end{equation*}
\(T(x) = x \in \mathbb{R}^k\), \(\eta(\theta)=\eta(p) = (ln(p_1), ..., ln(p_k))\), $H(\theta)=0$, $a(x_i)=\frac{N!}{x_{i,1}!...x_{i,k}!}$ and $q=k$. However, introducing \(\eta_k(\eta) = ln(p_k)\) poses a challenge, as it would render the model non-identifiable. This non-identifiability arises from \(H(\eta) = 0\) and \(\nabla H(\eta) = 0\), implying that the entire \(\mathbb{R}^k\) becomes a critical point of \(l_n(\eta)\) in Eq. \ref{ln eta}.

\section{Case studies}
\label{case studie}
To illustrate our formulas we'll use the Weibull distribution. This density is used in various applications, including image processing \cite{medical}, survival analysis \cite{survival}, and extreme value theory \cite{extreme value theory}. To obtain the Lehmer and Holder means as MWLE, we focus on the scenario of independent variables. In the following sections, we present the Lehmer case first and then proceed to the H\"older case.
\newline
Let be $X=(X_1,X_2,X_3) \sim Weibull(\lambda,k)$ with $\lambda=(\lambda_1,\lambda_2,\lambda_3)\in\mathbb{R}_{\geq 0}^3$ of shape $k=(k_1,k_2,k_3)\in\mathbb{R}_{\geq 0}^3$. As $X$ has independent random variable elements, so:
\begin{equation}
    f(x_i)=\prod_{j=1}^3 f(x_{i,j})=\prod_{j=1}^3(\frac{k_j}{\lambda_j}(\frac{x_{i,j}}{\lambda_j})^{k_j-1}e^{-\frac{x_{i,j}}{\lambda_j}})
\end{equation}
The Weibull PDF belongs to the exponential family when:
 \begin{equation}
    \begin{split}
        &a(x_i)=x_{i,1}^{k_1-1}x_{i,2}^{k_2-1}x_{i,3}^{k_3-1}
        \\& exp(H(\lambda))=\frac{k_1k_2k_3}{\lambda_1^{k_1}\lambda_2^{k_2}\lambda_3^{k_3}}
        \\&
        T(x_i)=\begin{pmatrix}
            x_{i,1}^{k_1}\ ,&x_{i,2}^{k_2}\ ,&x_{i,3}^{k_3}
        \end{pmatrix}
        \\& 
        \eta(\lambda)=\begin{pmatrix}
            \lambda_1^{-k_1}\ ,&\lambda_2^{-k_2}\ ,&\lambda_3^{-k_3}
        \end{pmatrix}
    \end{split}
    \label{Weibull expo}
\end{equation}
As a reminder, $\forall i\in\{1,2,3\}$ the moments of $X_i$ is given by:
\begin{equation}
    \forall t\geq 0, \  \mathbb{E}_{X_i\sim\eta_i}[X_i^t]=\lambda_i^t \Gamma(1+\frac{t}{k_i})
    \label{moments}
\end{equation}
Let us recall that this PDF belongs to the minimal exponential family.
 In this case, the vector $\eta$ is then unique (see Appendix \ref{minimal_expo}).
In order to validate the use of the H\"older and Lehmer mean families as MWLE, we use the data ${\cal D}$ of the U.S. Senate statewide
1976–2020 elections \cite{Harvard} to illustrate our calculations: The
dataset ${\cal D}$ is composed of three variables representing the proportion of votes for each political parties DEMOCRATS, OTHERS and REPUBLICANS. ”OTHERS” being the less important political parties regrouped together. The proportions
of votes a political party obtains are the random variables:
$X_1$ for DEMOCRATS, $X_2$ for REPUBLICANS and $X_3$ for
OTHER.
There were 24 elections between 1967-2020; that is according to Eq. \ref{Dataset weight}, the data ${\cal D}$ is a 24x3 matrix. In what follows, we estimate the average vote for each political party.

\subsection{Lehmer's mean}
Let be ${\cal D} = {\cal H} \times {\cal U}$
(see Eq. \ref{Dataset weight}). 
 Let us recall that the Lehmer mean family is derived when the shape parameter $k=1$. So, $T_j(x_{i,j})=x_{i,j}$,  
  $u(x_{i,j})=w(x_{i,j}) x_{i,j}^{\beta_j-1}$ with $w(x_{i,j})=1$. In this case Eq. \ref{CP} and Eq. \ref{Weibull expo} is rewritten as:
\begin{equation}
\begin{split}
   r(\eta)&=\begin{pmatrix}
    \frac{\sum_{i=1}^n u(x_i) x_{i,1}}{\sum_{i=1}^n u(x_i)}\ ,&\frac{\sum_{i=1}^n u(x_i) x_{i,2}}{\sum_{i=1}^n u(x_i)}\ ,&\frac{\sum_{i=1}^n u(x_i) x_{i,3}}{\sum_{i=1}^n u(x_i)}
        \end{pmatrix} 
        \\&=\begin{pmatrix}
    \mathbb{E}_{X_1\sim\eta_1}[X_1]\ ,&\mathbb{E}_{X_2\sim\eta_2}[X_2]\ ,&\mathbb{E}_{X_3\sim\eta_3}[X_3]
        \end{pmatrix}
        \end{split}
\end{equation}
Substituting $t=k_1=1$ in Eq. \ref{moments}, $r(\eta)$ is rewritten as:
\begin{equation}
            \begin{split}
                    r(\eta) =&\begin{pmatrix}  \lambda_1\ , &\lambda_2\ , &\lambda_3
        \end{pmatrix}
        \\=&\lambda
\end{split}
\label{lehmer estimation}
\end{equation}
Then:

\begin{equation*}
 \lambda\sim
  \begin{pmatrix}
     \frac{\sum_{i=1}^n x_{i,1}^{\beta_1}}{\sum_{i=1}^n x_{i,1}^{\beta_1-1}} \ , & \frac{\sum_{i=1}^n x_{i,2}^{\beta_2}}{\sum_{i=1}^n x_{i,2}^{\beta_2-1}} \ ,& \frac{\sum_{i=1}^n x_{i,3}^{\beta_3}}{\sum_{i=1}^n x_{i,3}^{\beta_3-1}}
  \end{pmatrix}   
 \end{equation*}
 In Fig. \ref{lehmer estimation representation}, on average, the Democrats ($\lambda_1$) are more likely
to be voted in the U.S Senate than the Republicans ($\lambda_2$)
or the Others ($\lambda_3$). The average here is not a specific value, but rather a function of $\beta_i$, ranging from the minimum of our data sample as $\beta_i$ approaches $-\infty$, to its maximum value as $\beta_i$ approaches $+\infty$. The choice of a value to express the tendency may depend on external factors, such as favoring certain elections over others. However, regardless of the value of $\beta_i$, there are invariants such as the order relation of the lambdas ($\lambda_1 > \lambda_2 > \lambda_3$) or that a lambda function is non-decreasing in $\beta$.

\subsection{H\"older's mean}
The H\"older's mean as MWLE can be found here because $T_j(x_{i})=x_{i,j}^{k_j}$,   $u(x_{i})=w(x_{i}) $.
We use Eq. \ref{CP} and \ref{Weibull expo} to find the first condition to have the MWLE:
\label{A section}
\begin{equation}
\begin{split}
   r(\eta)&=\begin{pmatrix}
    \frac{\sum_{i=1}^n u(x_i) x_{i,1}^{k_1}}{\sum_{i=1}^n u(x_i)}\ ,&\frac{\sum_{i=1}^n u(x_i) x_{i,2}^{k_2}}{\sum_{i=1}^n u(x_i)}\ ,&\frac{\sum_{i=1}^n u(x_i) x_{i,3}^{k_3}}{\sum_{i=1}^n u(x_i)}
        \end{pmatrix} 
        \\&=\begin{pmatrix}
    \mathbb{E}_{X_1\sim\eta_1}[X_1^{k_1}]\ ,&\mathbb{E}_{X_2\sim\eta_2}[X_2^{k_2}]\ ,&\mathbb{E}_{X_3\sim\eta_3}[X_3^{k_3}]
        \end{pmatrix}
\end{split}
\label{restimator}
\end{equation}
We notice that the components of this vector are non-centered moments of order $k_j$. Let us recall that the weighted H\"older mean family when $u(x_i)=w(x_i)=1$.  Eq. \ref{restimator} is rewritten as: 
\begin{equation}
\begin{split}
    r(\eta)=\begin{pmatrix}
    \frac{1}{n}\sum_{i=1}^n x_{i,1}^{{k_1}} \ ,&\frac{1}{n}\sum_{i=1}^n x_{i,2}^{k_2}\ ,&\frac{1}{n}\sum_{i=1}^n  x_{i,3}^{k_3} \end{pmatrix}
    \end{split}
    \end{equation}
    Substituting $t=k_j$ in Eq. \ref{moments}, $r(\eta)$ is rewritten as:
    \begin{equation}
        \begin{split}
r(\eta)= \begin{pmatrix}
    \lambda_1^{k_1} \ ,& \lambda_2^{k_2} \ ,&\lambda_3^{k_3}
    \end{pmatrix}
\label{eq.r}
\end{split}
\end{equation}
The estimators of the proportions to vote for $X_1$, $X_2$ and $X_3$ are then thanks to Eq. \ref{eq.r}:
\begin{equation}
    \lambda_1=(\frac{1}{n}\sum_{i=1}^n x_{i,1}^{{k_1}})^{\frac{1}{k_1}},\ \lambda_2=(\frac{1}{n}\sum_{i=1}^n x_{i,2}^{k_2})^{\frac{1}{k_2}}\ \text{and}\ \lambda_3=(\frac{1}{n}\sum_{i=1}^n  x_{i,3}^{k_3})^{\frac{1}{k_3}}
    \label{Horder estimate 1}
\end{equation}
The figure \ref{fig Holder} represents the three estimators of equation \ref{Horder estimate 1}. It shows that $\lambda_1>\lambda_2>\lambda_3$, so the prominent political parties in the Senate is on average the DEMOCRATS, followed by the REPUBLICANS and then the OTHERS. Here, the weight is $u(x_i)=1$, the data have the same weight. But the h\"older's mean parameter $k$ is also the shape parameter of the Weibull distribution, so it can adjust to our data. Besides, the MWLE which depends on $k_i\in\mathbb{R}_{\geq0}$ is an increasing function. Furthermore, the H\"older's mean goes from the harmonic mean of our data sample when $k_i=0$ and tends to the maximum of our data sample when $k_i$ tends to $+\infty$.
\section{Comparison}
\label{comparaison}
We have shown that the H\"older mean and the Lehmer mean are MWLEs. There are significant differences between them. In the H\"older mean, the parameter $k_i$ is always positive, while in the case of the Lehmer mean, the parameter $\beta_i$ can be negative. Consequently, depending on $k_i$, the H\"older mean is located between the harmonic mean and the maximum observation. In contrast, the Lehmer mean is located between the minimum observation and the maximum observation. The parameters of these two estimators fulfill different roles: $k_i$ is a shape parameter of the PDF, calculable from the data, but $\beta$ implements a data selection. Both estimators are non-decreasing in $k$ for H\"older and $\beta$ for Lehmer. However, the H\"older mean increases at a slower rate than the Lehmer mean over the entire range from 0 to $+\infty$, as indicated by their slopes (see Fig. \ref{lehmer estimation representation} and Fig. \ref{fig Holder}).

\section{Conclusions}
Among the infinite definitions of the mean, there are the families of Lehmer and H\"older means. These mean families are used quite often. Recently, their link with maximum likelihood estimators has been established only in the case of univariate probability densities. In this article, we have shown that these two families of means are maximum likelihood estimators in the case of weighted data and PDFs belonging to the multivariate exponential family. We therefore offer a probabilistic interpretation of these mean families in multidimensional spaces. 

\appendix
\section{The MWLE}
\label{calcul_det}
The objective of this appendix is to derive the MWLE.
Let be ${\cal D} = {\cal H}\times{\cal U} $ a weighted dataset of n observation sampled from a random vector $X=(X_1,...,X_k)$ of joint distribution $\mathbb{P}(\theta)$ of parameter $\theta \in \Theta \subset \mathbb{R}^k$ (see Eq. \ref{Dataset weight}).
We define the weighted log-likelihood function relatively to the natural parameter $\eta$ is given by:
\begin{equation*}
     l_n(\eta)=\sum_{i=1}^n u(x_i)[\ln(a(x_i)) +<\eta,T(x_i)>  - H(\eta))]
\end{equation*}
$l_n()$ is a $C^2$ function relatively to $\eta$ because it is the log of the exponential function. The parameters $\eta$, $T()$ and $H()$ are defined in Eq. \ref{pdfgexp}. We will now maximize the weighted log-likelihood. To this end, when differentiating partially by $\eta_j$, $j \in \{1,...,q\}$ we find:
\begin{equation}
       \frac{\partial l_n}{\partial \eta_j}(\eta) = \sum_{i=1}^n u(x_i)(<\frac{\partial \eta}{\partial \eta_j},T(x_i)> - \frac{\partial H}{\partial \eta_j}(\eta)))
   \label{ln deriv}
\end{equation}
Where:
\begin{equation*}
    \frac{\partial H}{\partial \eta_j}(\eta) = \frac{\int_{\mathbb{R}_{\geqslant0}^k}<\frac{\partial \eta}{\partial \eta_j},T(x)>a(x)exp(<\eta,T(x)>)dx}{\int_{\mathbb{R}_{\geqslant0}^k}a(x)exp(<\eta,T(x)>)dx}
\end{equation*}
\begin{equation*}
 = \frac{\int_{\mathbb{R}_{\geqslant0}^k}T_j(x)a(x)exp(<\eta,T(x)>)dx}{\int_{\mathbb{R}_{\geqslant0}^k}a(x)exp(<\eta,T(x)>)dx}
\end{equation*}
\begin{equation}
        = \mathbb{E}_{X\sim\eta}[T_j(X)]
\end{equation}
Consequently, we can write:
\begin{equation}
    \grad H(\eta)=\mathbb{E}_{X\sim\eta}[T(X)]
\end{equation}
\begin{equation}
\begin{split}
    \grad l_n(\eta)= \sum_{i=1}^n u(x_i)(T(x_i) - \mathbb{E}_{X\sim\eta}[T(X)])
\end{split} 
\end{equation}
The first order optimality condition leads to the critical point $\eta*$:
\begin{equation}
   \begin{split}
       r(\eta*) &= \mathbb{E}_{X\sim\eta*}[T(X)] \\&= \frac{\sum_{i=1}^n u(x_i)T(x_i)}{\sum_{i=1}^n u(x_i)} 
   \end{split}   
\end{equation} 
Let us now, examine the second order optimality condition. For this end, let's compute the Hessian matrix $\forall \eta \in \mathbb{R}^q$. $\forall (i,j)\in\{1,...,q\}$, as $l_n()$ is a $C^2$ function, by Cauchy-Schwartz:
\begin{equation}
   \begin{split}
       \frac{\partial^2 l_n}{\partial\eta_i \partial\eta_j}(\eta) &= \frac{\partial^2 l_n}{\partial\eta_j \partial\eta_i}(\eta) \\ & =-(\sum_{i=1}^n u(x_i))\frac{\partial \mathbb{E}_{X\sim\eta}[T_i(X)]}{\partial \eta_j}(\eta)\\ & = -(\sum_{i=1}^n u(x_i))\frac{\partial \mathbb{E}_{X\sim\eta}[T_j(X)]}{\partial \eta_i}(\eta)
   \end{split} 
\end{equation} 
\begin{equation*}
\begin{split}
     = &-(\sum_{i=1}^n u(x_i))[\int_{\mathbb{R}_{\geqslant0}^k}T_i(x) T_j(x) a(x)exp(<\eta,T(x)>-H(\eta))  dx \\&
- (\int_{\mathbb{R}_{\geqslant0}^k}T_i(x) a(x)exp(<\eta,T(x)>-H(\eta))dx)\\&(\int_{\mathbb{R}_{\geqslant0}^k} T_j(x) a(x)exp(<\eta,T(x)>-H(\eta))dx)]
\end{split}
\end{equation*}
\begin{equation*}
    \begin{split}
&=-(\sum_{i=1}^n u(x_i))[\mathbb{E}_{X\sim\eta}[T_i(X)T_j(X)] - \mathbb{E}_{X\sim\eta}[T_i(X)]\mathbb{E}_{X\sim\eta}[T_j(X)]] \\&=-(\sum_{i=1}^n u(x_i))[Cov_{X\sim\eta}(T_i(X),T_j(X))]
    \end{split}
\end{equation*}
So,
\begin{equation}
    Hess_{l_n} (\eta)=-(\sum_{i=1}^n u(x_i))K_{T(X),T(X)}
    \label{Condition second ordre}
\end{equation}
Where $K_{T(X),T(X)}$ is the covariance matrix of $T(X)$ and $\sum_{i=1}^n u(x_i)$ the total weight which is positive and non-null. A covariance matrix is always a semi-definite positive matrix. So the Hessian matrix in Eq. \ref{Condition second ordre} is a negative semi-definite Matrix for any $\eta$. It follows that any critical point is a maximum. There are only two possible forms of the Hessian: 1) positive definite, in this case ln() has a single maximum; 2) degenerates, in this case it is possible that a maximum degenerates (i.e., an infinite  number of maximums, a sadle point). In appendix~\ref{minimal_expo}, we will identify the class of the exponential family which has a single maximum. If so, there is a one-to-one mapping between r() and $\eta*$. The existence of unique maximum $\theta$* requires that $\eta$() is bijective on the space parameter $\Theta$.
\begin{equation*}
    \theta\sim\eta^{-1}(r^{-1}(\frac{\sum_{i=1}^n u(x_i)T(x_i)}{\sum_{i=1}^n u(x_i)}))
\end{equation*}
With $r^{-1}()$ and $\eta^{-1}()$ representing the inverse of the components of the vectors $r()$ and $\eta()$.

\section{Estimator uniqueness}
\label{minimal_expo}
The objectif of this appendix is to find the sub-class of PDFs such that the weighted likelihood $l_n()$ have an unique maximum.  \newline
\underline{Theorem:} Let be $\cal{D}=\cal{H}\times\cal{U}$ (see Eq. \ref{Dataset weight}) and $X=(X_1,...,X_k)$ a random vector following a distribution $\mathbb{P}(\eta)$ with $\eta \in \mathbb{R}^q$, whose PDF belongs to a minimal exponential family, such as $f(x|\eta)=a(x)\exp(<\eta, T(x)>-H(\eta))$.
If a solution of the equation, 
\begin{equation}
    r(\eta)= \frac{\sum_{i=1}^n u(x_i)T(x_i)}{\sum_{i=1}^n u(x_i)}  \label{Extremum point}
\end{equation}
 exists then the weighted log-likelihood density function of X has a unique maximum, and $l_n()$ becomes bijective at $\eta*$.
 \newline
\underline{Proof:}
\newline
Let us recall that the covariance matrix is semi-definite positive; that is  $\forall y\in\mathbb{R}^q\setminus\{0\}, \ y^tK_{T(X),T(X)}y>0$. We would like find the class of exponential family PDFs leading to  a strict inequality. In other words, that case $y^tK_{T(X),T(X)} y=0$ cannot happens for $y \neq 0$.
This can be demonstrated by proving that $\forall y\in\mathbb{R}^q\setminus\{0\}, \  y^tK_{T(X),T(X)}y > 0$. Let be $y=(y_1,...,y_q)\in\mathbb{R}^q$ and $\eta\in\mathbb{R}^q$,
\begin{equation*}
y^tK_{T(X),T(X)}y 
\end{equation*}
 \begin{equation*}
    =  \begin{pmatrix}
        y_1 \\ \vdots \\ y_q
    \end{pmatrix}^t \begin{pmatrix}
       \sum_{j=1}^q y_j Cov_{X\sim\eta}(T_1(X),T_j(X)) \\ \vdots\\ \sum_{j=1}^q y_j Cov_{X\sim\eta}(T_q(X),T_j(X))
    \end{pmatrix} 
 \end{equation*}
\begin{equation*}
= \begin{pmatrix}
        y_1 \\ \vdots \\ y_q
    \end{pmatrix}^t \begin{pmatrix}
        Cov_{X\sim\eta}(T_1(X),\sum_{j=1}^q y_j T_j(X)) \\ \vdots\\  Cov_{X\sim\eta}(T_q(X),\sum_{j=1}^q y_jT_j(X))
    \end{pmatrix} 
\end{equation*}
\begin{equation*}
    =   \sum_{i=1}^q y_i Cov_{X\sim\eta}(T_i(X),\sum_{j=1}^q y_j T_j(X)) 
    \end{equation*}
    \begin{equation*}
=Cov_{X\sim\eta}(\sum_{i=1}^q y_iT_i(X),\sum_{j=1}^q y_j T_j(X)) 
\end{equation*}

\begin{equation*}
  =  Var_{X\sim\eta}(\sum_{j=1}^p y_j T_j(X)) \label{vartx}
\end{equation*}
$Var_{X\sim\eta}(\sum_{j=1}^p y_j T_j(X))=0$ implies that $\sum_{j=1}^p y_j T_j(X)$ is constant. Let's consider the minimal exponential family, a subclass of the exponential family. In this family, we cannot have coefficients $y\neq0$ such that $\sum_{j=1}^p y_j T_j(X)$ is constant. Therefore, for this subclass, $ y^tK_{T(X),T(X)}y >0$ for $y \neq 0$. It follows that the Hessian is negatively defined (i.e. $l_n()$ is strictly convex), and thus the critical point is the unique maximum. 


\begin{figure} 
\begin{center}
        \centering
        \includegraphics[width=7cm]{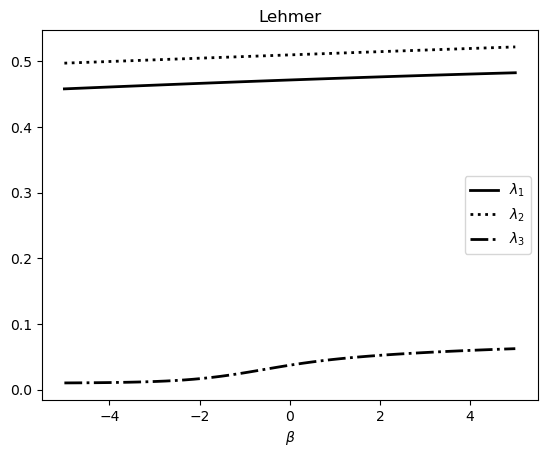}       
\caption{MWLE of $\lambda$ as the Lehmer's mean as function of $\beta$}
        \label{lehmer estimation representation}
\end{center}
\end{figure}

\begin{figure} 
\begin{center}
\includegraphics[width=7cm]{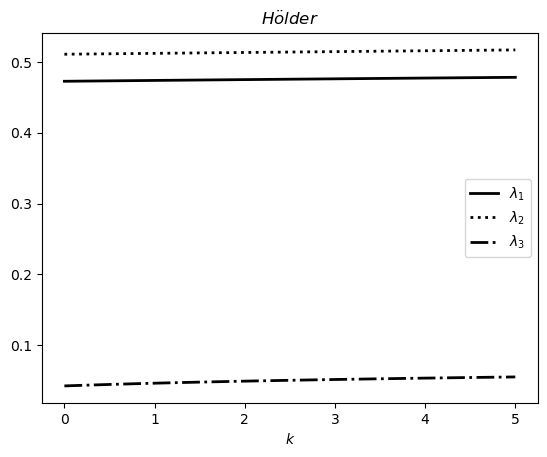}
\caption{MWLE of $\lambda$ as the H\"older's mean as function of $k$}
\label{fig Holder}
\end{center}
\end{figure}



\bibliographystyle{elsarticle-num}

\end{document}